\documentclass[11pt]{article}
\usepackage{graphics}
\usepackage{latexsym}
\newtheorem{theorem}{Theorem}

\newtheorem{proposition}{Proposition}
\newtheorem{problem}{Problem}

\newtheorem{lemma}[proposition]{Lemma}

\newtheorem{corollary}{Corollary}
\newtheorem{corollaryth}{Corollary}

\newcommand{\qed}{\ \hfill\mbox{$\Box$}\vspace{\baselineskip}}
\newenvironment{proof}{\noindent {\bf Proof:}}{{\qed}}

\begin{document}

\title{On Gale and braxial polytopes}

\author{Margaret M. Bayer\thanks{Supported in part by a grant from the
        University of Kansas General Research Fund}\\
        Department of Mathematics\\
        University of Kansas\\
        Lawrence KS 66045-7523 USA\\ \and
        Tibor Bisztriczky\thanks{Supported in part by a Natural Sciences and
        Engineering Research Council of Canada Discovery Grant}\\
        Department of Mathematics and Statistics\\
        University of Calgary\\
        Calgary, Alberta, T2N 1N4 Canada}

\date{September 2006}
\maketitle

\begin{abstract}
Cyclic polytopes are characterized as simplicial polytopes satisfying Gale's
evenness condition (a combinatorial condition on facets relative to a fixed
ordering of the vertices).  Periodically-cyclic polytopes are polytopes
for which certain subpolytopes are cyclic.  Bisztriczky discovered a class
of periodically-cyclic polytopes that also satisfy Gale's evenness condition.
The faces of these polytopes are braxtopes, a certain class of nonsimplicial
polytopes studied by the authors.  In this paper we prove that the
periodically-cyclic Gale polytopes of Bisztriczky are exactly the polytopes
that satisfy Gale's evenness condition and are braxial (all faces are
braxtopes).  The existence of other periodically-cyclic Gale polytopes
is open.

\end{abstract}

\section{Introduction}
We recall that cyclic polytopes have a totally ordered set of vertices (vertex 
array) that satisfies Gale's Evenness Condition and yields a complete 
description of their facet structure.
One seeks to generalize this class of polytopes due to their
important combinatorial properties and to their many applications in various 
branches of mathematics and science.
Of specific significance are generalizations that are nonsimplicial and that 
exhibit constructions other than products, pyramids, prisms and so forth.
Bicyclic 4-polytopes, ordinary $d$-polytopes and certain periodically-cyclic 
Gale $d$-polytopes are examples of such generalizations.
It is noteworthy that these are also polytopes with explicit facet structures.

Our present interest is characterizations of these polytopes that do not invoke 
their constructions or facet structures.
For example: cyclic polytopes may be characterized as Gale and simplicial 
polytopes or as neighbourly polytopes with the same number of universal edges as
vertices or as polytopes that have only cyclic subpolytopes.
With the observation that multiplices are generalizations of simplices, ordinary
polytopes may be characterized as Gale and multiplicial polytopes.
With the knowledge that braxtopes are also generalizations of simplices and that
they were discovered as facets of certain periodically-cyclic Gale polytopes, it
is natural to ask if there is a characterization of periodically-cyclic Gale 
polytopes as polytopes that are Gale and braxial.
In the following, we determine that all Gale and braxial $d$-polytopes are 
periodically-cyclic for $d\ge 5$.

\section{Definitions and background}
Let $Y$ be a set of points in ${\bf R}^d$, $d\ge 1$.
Then $[Y]$ and $\langle Y\rangle$ denote, respectively, the convex hull and the affine 
hull of $Y$.
If $Y=\{y_1,y_2,\ldots, y_s\}$ is finite, we set $[y_1,y_2,\ldots, y_s]=[Y]$
and $\langle y_1,y_2,\ldots, y_s\rangle =\langle Y\rangle$

Let $X=\{x_0,x_1,\ldots, x_n\}$ be a totally ordered set of $n+1$ points in 
${\bf R}^d$ with $x_i<x_j$ if and only if $i<j$.
We say that $x_i$ and $x_{i+1}$ are {\em successive} points, and if
$x_i<x_j<x_k$ then $x_j$ {\em separates}, or is {\em between}, $x_i$ and $x_k$.
Let $Y\subset X$.
Then $Y$ is a {\em Gale} subset of $X$ if any two points of $X\setminus Y$ are 
separated by an even number of points of $Y$.
Finally, $X$ is a {\em paired set} if it is the union of mutually disjoint
$\{x_i,x_{i+1}\}$.
As a rule, $S_m$ denotes a paired set of $m$ points with $S_0$ denoting the
empty set.

Let $P\subset {\bf R}^d$ be a (convex) $d$-polytope.
For $-1\le i\le d$, let ${\cal F}_i(P)$ denote the set of $i$-dimensional 
faces of $P$ with $f_i(P)=|{\cal F}_i(P)|$.
For convenience, let
$$\mbox{${\cal V}(P)={\cal F}_0(P)$, ${\cal E}(P)={\cal F}_1(P)$, and
${\cal F}(P)={\cal F}_{d-1}(P)$.}$$
We assume familiarity with the basic definitions and concepts concerning
polytopes (see \cite{grunbaum,ziegler}), and we cite two results necessary
for our presentation from \cite{grunbaum} and \cite{klee-minty}, respectively.

\begin{lemma}\label{inductive}
Let $P'$ and $P$ be $d$-polytopes in ${\bf R}^d$ such that $P=[P',x]$ for some
point $x\in{\bf R}^d\setminus P'$.
Let $G$ be a face of $P'$ and 
${\cal F}(G,P')=\{F\in{\cal F}(P')\,|\:G\subseteq F\}$.
Then
\begin{enumerate}
\item $G$ is a face of $P$ if and only if $x$ is beneath some 
      $F\in{\cal F}(G,P')$;\label{1a}
\item $[G,x]$ is a face of $P$ if and only if either 
      $x\in\langle G\rangle$ or $x$ is
      beneath $F'$ and beyond $F''$ for some $F'$ and $F''$ in
      ${\cal F}(G,P')$; and\label{1b}
\item if $G$ is a face of $P$ and $[G,x]$ is not a face of $P$, then $x$ is
      not beyond any $F\in{\cal F}(G,P')$.\label{1c}
\end{enumerate}
\end{lemma}

\begin{lemma}
If the facet system of one $d$-polytope is contained in the facet system of 
another $d$-polytope, then the two $d$-polytopes are combinatorially 
equivalent.
\end{lemma}

Let ${\cal V}(P)=\{x_0,x_1,\ldots, x_n\}$, $n\ge d$.
We let $x_i<x_j$ if and only if $i<j$, and call $x_0<x_1<\cdots <x_n$ a 
{\em vertex array} of $P$.
Let $G\in{\cal F}_i(P)$, $1\le i\le d-1$, and
$${\cal V}(G)=G\cap {\cal V}(P)=\{y_0,y_1,\ldots, y_m\};$$
that is, each $y_s$ is some $x_t$.
Then $y_0<y_1<\cdots <y_m$ is the vertex array of $G$ if it is the ordering
induced by $x_0<x_1<\cdots<x_n$.
Finally, $P$ with $x_0<x_1<\cdots<x_n$ is {\em Gale} (with respect to the
vertex array) if ${\cal V}(F)$ is a Gale set for each $F\in {\cal F}(P)$.
We recall from \cite{gale,grunbaum} that $P$ is {\em cyclic} if it is
simplicial, and Gale with respect to some vertex array.

From \cite{bisz-mult}, $P$ is a {\em multiplex} if there is a vertex array,
say, $x_0<x_1<\cdots<x_n$ such that 
$${\cal F}(P)=\{[x_{i-d+1},\ldots, x_{i-1},x_{i+1},\ldots, x_{i+d-1}]\,|\:
0\le i\le n\}$$
under the convention: $x_t=x_0$ for $t\le 0$ and $x_t=x_n$ for $t\ge n$.
A $d$-multiplex is a natural generalization of a $d$-simplex.
Next, $P$ is {\em multiplicial} if each facet of $P$ is a $(d-1)$-multiplex
with respect to the ordering induced by a fixed vertex array of $P$.
Finally, $P$ is {\em ordinary} if it is Gale and multiplicial with 
respect to some vertex array.
We note from \cite{bisz-ord} that if $P$ is an ordinary $d$-polytope
and $d\ge 4$ then there is a complete description of ${\cal F}(P)$;
furthermore, if $d$ is even then $P$ is cyclic.

\pagebreak

From \cite{bisz-pcg}, $P$ is {\em periodically-cyclic} if there is a vertex 
array, say, $x_0<x_1<\cdots<x_n$ and an integer $k$,
$d+2\le k\le n+1$, such that 
\begin{itemize}
\item $[x_{i+1},x_{i+2},\ldots, x_{i+k}]$ is a cyclic $d$-polytope with the
      induced vertex array for $-1\le i\le n-k$, and
\item $[x_{i+1},x_{i+2},\ldots, x_{i+k}, x_{i+k+1}]$ is not cyclic for
      $-1\le i\le n-k-1$.
\end{itemize}
The integer $k$ is the {\em period} of $P$.
We note that if $k=n+1$, then $P$ is cyclic.

Let $P$ be a periodically-cyclic $d$-polytope with $x_0<x_1<\cdots<x_n$.
Since the condition that 
$[x_{i+1},x_{i+2},\ldots, x_{i+k+1}]$ is not a cyclic $d$-polytope may be
satisfied in numerous ways, it follows that $P$ may be one of many 
combinatorial types.
This observation remains valid even under the added assumption that $P$ is Gale
with $x_0<x_1<\cdots<x_n$; see, for example, the bicyclic 4-polytopes in 
\cite{smilansky} that are Gale and periodically-cyclic \cite{bisz-pcg4}.
For $d>4$, we know at present one class of realizable periodically-cyclic
$d$-polytopes \cite{bisz-pcg}.

\begin{proposition} \label{cyclic}
Let $d\ge 6$ be even and $k\ge d+2$.
Let $P_{k-1}$ be a cyclic $d$-polytope in ${\bf R}^d$ with
vertex array $x_0<x_1<\cdots<x_{k-1}$.
Then there exist a sequence of points $x_n\in{\bf R}^d$ 
and a sequence of polytopes $P_n=[P_{n-1},x_n]$ ($n\ge k$) such that 
\begin{itemize}
\item $x_n\in\langle x_0,x_{n-k+1},x_{n-k+2},x_{n-1}\rangle$;
\item $x_n$ is beyond each $F\in{\cal F}(P_{n-1})$ with the property that 
      $F\cap[x_0,x_{n-k+1},x_{n-1}]=[x_0,x_{n-1}]$; and 
\item $x_n$ is beneath each $F\in{\cal F}(P_{n-1})$ with the property that
      $[x_0,x_{n-k+1},x_{n-k+2},x_{n-1}]\not\in F$ and
      $F\cap[x_0,x_{n-k+1},x_{n-1}]\ne[x_0,x_{n-1}]$.
\end{itemize}
Each polytope constructed in this manner is Gale and periodically-cyclic with
respect to $x_0<x_1<\cdots<x_n$ and with period $k$.
\end{proposition}

In \cite{bisz-pcg} there is an explicit combinatorial description,
depending only on $k$ and $n$, of the facets of these polytopes.
It is tedious but straightforward to check that these facets are of the 
following combinatorial type.

\pagebreak

\noindent {\bf Definition.}
Let $Q\subset{\bf R}^e$ be an $e$-polytope with 
${\cal V}(Q)=\{y_0,y_1,\ldots, y_m\}$, $m\ge e\ge 3$.
Then $Q$ is an {\em $e$-braxtope} if there is a vertex array, say,
$y_0<y_1<\cdots<y_m$ such that 
$${\cal F}(Q)=\{T_0,T_1,\ldots, T_{m-e+1},E_2,E_3,\ldots, E_m\}$$
with 
$$\mbox{$T_i=[y_i,y_{i+1},\ldots, y_{i+e-1}]$ for $0\le i\le m-e+1$},$$
and
$$\mbox{$E_j=[y_0,y_{j-e+2},\ldots, y_{j-1},y_{j+1},\ldots, y_{j+e-2}]$
for $2\le j\le m$},$$
under the convention that 
$y_t=y_0$ for $t\le 0$ and $y_t=y_m$ for $t\ge m$.
For the sake of completeness, an $e$-braxtope is an $e$-simplex for
$0\le e\le 2$.

\vspace{1\baselineskip}

It is clear that an $e$-braxtope with $e+1$ vertices is an $e$-simplex,
and that braxtopes are generalizations of simplices.

Finally, a $d$-polytope $P$ is {\em braxial} if each facet of $P$ is a 
$(d-1)$-braxtope with respect to the ordering induced by a fixed vertex
array of $P$.
As noted above, the periodically-cyclic $d$-polytopes constructed via 
Proposition~\ref{cyclic} are Gale and braxial.
The following result from \cite{bayer-bisz} enables us to prove that Gale
and braxial $d$-polytopes are periodically-cyclic for $d\ge 5$.

\begin{lemma}\label{edges}
Let $Q$ be an $e$-braxtope with $y_0<y_1<\cdots<y_m$, $m\ge e\ge 3$.
Then with the convention that 
$y_t=y_0$ for $t\le 0$ and $y_t=y_m$ for $t\ge m$:
\begin{enumerate}
\item $[y_0,y_t]\in{\cal E}(Q)$ for $1\le t\le m$. \label{4a}
\item $[y_1,y_t]\in{\cal E}(Q)$ if and only if $t\in\{0,2,3,\ldots, e\}$.
      \label{4b}
\item $[y_t,y_m]\in{\cal E}(Q)$ if and only if $t\in\{0,m-e+1,\ldots, m-1\}$.
      \label{4c}
\item For $2\le s\le m-1$, $[y_s,y_t]\in{\cal E}(Q)$ if and only if 
      $t\in\{0,s-e+1,\ldots, s-1,s+1,\ldots, s+e-1\}$.
      \label{4d}
\item $[y_0,y_t,y_{t+1},y_{t+e-1},y_{t+e}]\in{\cal F}_3(Q)$ for $1\le t\le m-e$.
      \label{4e}
\item $\{y_t,y_{t+1},\ldots, y_{t+e}]$ is an affinely independent set for
      $0\le t\le m-e$.
      \label{4f}
\end{enumerate}
In addition, $Q$ is braxial and if $m\ge e+1$ then 
$[y_0,y_1,\ldots, y_{m-1}]$ is an $e$-braxtope with 
$y_0<y_1<\cdots<y_{m-1}$.
\end{lemma}

\section{Gale and braxial polytopes}
Henceforth, we assume that $P$ is a Gale and braxial $d$-polytope with
respect to $x_0<x_1<\cdots<x_n$, $n\ge d+1$ and $d\ge 3$.
We simplify our notation by assuming also that $F$ and $F'$ always denote
facets of $P$ with the following properties (see Section~2 with $e=d-1$):

\vspace{-\baselineskip}

\begin{eqnarray}\label{ridges1}
& & \mbox{$y_0<y_1<\cdots<y_m$ is the induced vertex array of $F$ and}\\
& & {\cal F}_{d-2}(F)= \{T_0,T_1,\ldots, T_{m-d+2},E_2,E_3,\ldots, E_m\},\;
m\ge d-1.\nonumber
\end{eqnarray}

\vspace{-\baselineskip}

\begin{eqnarray}\label{ridges2}
& & \mbox{$z_0<z_1<\cdots<z_u$ is the induced vertex array of $F'$ and}\\
& & {\cal F}_{d-2}(F')= \{T'_0,T'_1,\ldots, T'_{u-d+2},E'_2,E'_3,\ldots, E'_u\},
\; u\ge d-1.\nonumber
\end{eqnarray}

In the next two proofs we distinguish the facets of a braxtope by the 
number of vertices.
In particular, $E_3$ and $E_{m-2}$ are the only facets of the $(d-1)$-braxtope
$F$ having exactly $d$ vertices.

\begin{lemma}\label{simplex}
Let $F\in{\cal F}(P)$ with $x_0\not\in F$.
Then $F$ is a $(d-1)$-simplex.
\end{lemma}
\begin{proof}
The statement is trivial for dimension three, since all 2-dimensional
braxtopes are simplices.
Assume now that $d\ge 5$.
With reference to (\ref{ridges1}), we suppose that $m\ge d$ and seek a 
contradiction.

We note that $x_0\not\in F$ and the Gale property yield that 
$(y_0,y_1,y_2,y_3)=(x_{r-1},x_r,x_s,x_{s+1})$ for some $2\le r<s\le n-1$.
We consider 
$$E_3=[y_0,y_1,y_2,y_4,\ldots, y_d]=[x_{r-1},x_r,x_s,y_4,\ldots, y_d]$$
and $F'\in{\cal F}(P)$ with $z_0<z_1<\cdots<z_u$ such that $E_3=F\cap F'$.
Since $f_0(E_3)=d$, it follows that $E_3\in\{E'_3,E'_{u-2}\}$.
If $E_3=E'_3$ then $x_{r-1}<x_r<x_s<z_3<y_4<\cdots <y_d<\cdots$ is the vertex
array of $F'$, and $F'\cap\{x_0,x_{s+1}\}=\emptyset$.
Since $x_0$ and $x_{s+1}$ are separated by exactly three vertices of $F'$,
we have a contradiction of the Gale property.

Let $E_3=E'_{u-2}=[z_0,z_{u-d+1},\ldots, z_{u-3},z_{u-1},z_u]$.
We note that $z_0=y_0=x_{r-1}$, and hence, $z_{u-d+1}=y_1=x_r=z_1$ and
$x_{r-1}<x_r<x_s<y_4<\cdots<y_{d-2}<z_{u-2}<y_{d-1}<y_d$ is the vertex array
of $F'$.
Again, $x_0$ and $x_{s+1}$ are separated by exactly three vertices of $F'$.

The dimension four case is handled in a similar way.
\end{proof}

\begin{theorem}\label{theorema}
Let $P$ be a Gale and braxial $d$-polytope with respect to 
$x_0<x_1<\cdots<x_n$.
Let $d\ge 3$ be odd.
Then $P$ is a cyclic $d$-polytope with respect to $x_0<x_1<\cdots<x_n$.
\end{theorem}
\begin{proof}
We note that it is sufficient to prove that $P$ is simplicial.
This holds for dimension three because all 2-dimensional braxtopes are 
simplices.
So assume $d\ge 5$, and let $F\in{\cal F}(P)$ with $y_0<y_1<\cdots<y_m$.
We suppose that $m\ge d$, and seek a contradiction.

Since $m\ge d$, it follows from Lemma~\ref{edges} that $x_0=y_0$.
If $x_n=y_m$ then $T_1=[y_1,\ldots, y_{d-1}]=F\cap \tilde{F}$ for some
$\tilde{F}\in{\cal F}(P)$, $\tilde{F}\cap\{x_0,x_n\}=\emptyset$ and
$\tilde{F}$ is a $(d-1)$-simplex by Lemma~\ref{simplex}.
Since $d$ is odd and $x_0$ and $x_n$ are separated by the $d$ vertices of $F'$,
it follows that $x_n\ne y_m$.

Since $x_n\not\in F'$, we obtain from the Gale property that 
\begin{equation} \label{paired}
\mbox{$\{y_{m-d},y_{m-d+1},\ldots,y_{m-1},y_m\}$ is a paired set}
\end{equation}
and $\{y_{m-3},y_{m-2},y_{m-1},y_m\}=\{x_{r-1},x_r,x_s,x_{s+1}\}$ for
some $2\le r<s\le n-2$.
We consider
\begin{eqnarray*}\label{em-2}
E_{m-2}&=&[y_0,y_{m-d+1},\ldots, y_{m-3},y_{m-1},y_m]\\
       &=&[x_0,y_{m-d+1},\ldots, y_{m-4},x_{r-1},x_s,x_{s+1}]
\end{eqnarray*}
and $F'\in{\cal F}(P)$ with $z_0<z_1<\cdots<z_u$ such that $E_{m-2}=F\cap F'$.
Since $f_0(E_{m-2})=d$, it follows that $E_{m-2}\in\{E'_3,E'_{u-2}\}$.
If $E_{m-2}=E'_3$ then $z_0=y_0=x_0$, $z_1=y_{m-d+1}$ and 
$x_0<y_{m-d+1}<y_{m-d+2}<z_3<y_{m-d+3}<\cdots<
y_{m-4}<x_{r-1}<x_s<x_{s+1}<\cdots$ is the vertex array of $F'$.
If $E_{m-2}=E'_{u-2}$ then 
$z_0=y_0<\cdots<y_{m-d+1}<\cdots<y_{m-4}<x_{r-1}<z_{u-2}<x_s<x_{s+1}$ is the
vertex array of $F'$.
In case of the former, we note that $y_{m-d+2}$ and $y_{m-d+3}$ are not
successive vertices, a contradiction by (\ref{paired}).
In case of the latter, $F'\cap\{x_r,x_n\}=\emptyset$ and $x_r$ and $x_n$
are separated by exactly three vertices of $F'$.
\end{proof}

In view of Theorem~\ref{theorema}, we may now assume that $d\ge 4$ and even.
Our first task is to determine ${\cal E}(P)$, and we let
$${\cal V}_0={\cal V}_0(P)=\{x_j\in{\cal V}(P)\,|\:[x_0,x_j]\in{\cal E}(P)\}$$
and
$${\cal V}_i={\cal V}_i(P)=\{x_j\in{\cal V}(P)\,|\:\mbox{$x_j\ne x_0$ and
$[x_i,x_j]\in{\cal E}(P)$}\}, \; 1\le i\le n.$$

\begin{lemma}\label{V0}
${\cal V}_0(P)={\cal V}(P)\setminus\{x_0\}$.
\end{lemma}
\begin{proof}
We note that $x_j\in {\cal V}_0$ for some $2\le j\le n-1$, and that it is
sufficient to show that $\{x_{j-1},x_{j+1}\}\subset{\cal V}_0$.

Since $[x_0,x_j]\in{\cal E}(P)$, there are $F^*$ and $\tilde{F}$ in
${\cal F}(P)$ such that $[x_0,x_j]\subseteq F^*\cap\tilde{F}$,
$x_{j-1}\not\in F^*$ and $x_{j+1}\not\in\tilde{F}$.
Then $x_{j+1}\in F^*$ and $x_{j-1}\in \tilde{F}$ by the Gale property, and
$\{x_{j-1},x_{j+1}\}\subset {\cal V}_0$ by Lemma~\ref{edges}(\ref{4a}).
\end{proof}

\begin{lemma}\label{between}
Let $1\le p<q<r\le n$ and $[x_p,x_r]\in{\cal E}(P)$.
Then $x_q\in {\cal V}_p(P)\cap{\cal V}_r(P)$.
\end{lemma}
\begin{proof}
We note that it is sufficient to prove that 
$\{[x_p,x_{r-1}],[x_{p+1},x_r]\}\subset {\cal E}(P)$.

Since $[x_p,x_r]\in{\cal E}(P)$ and $x_p\ne x_0$, there is an $F\in{\cal F}(P)$
such that $[x_p,x_r]\subset F$ and $x_{p-1}\not\in F$.
Then $x_{p+1}\in F$ by the Gale property with, say, $x_r=y_s$ and 
$(x_p,x_{p+1})=(y_t,y_{t+1})$; see (\ref{ridges1}).
Now either $F$ is a simplex and $[x_{p+1},x_r]\in{\cal E}(P)$ or  $x_0=y_0$ by
Lemma~\ref{simplex}.
Let $x_0=y_0$.
Then $y_0\ne y_t$ and we apply Lemma~\ref{edges}(\ref{4c},\ref{4d}).
Specifically, $[y_t,y_s]=[x_p,x_r]\in{\cal E}(P)$ implies that 
$s-d+2\le t\le s-2$, whence $[x_{p+1},x_r]=[y_{t+1},y_s]\in{\cal E}(P)$.

In the case that $x_r\ne x_n$, a similar argument yields that
$[x_p,x_{r-1}]\in{\cal E}(P)$.

Let $[x_p,x_n]\in{\cal E}(P)$.
We note that if $[x_\ell,x_n]\in{\cal E}(P)$ implies that 
$[x_\ell,x_{n-1}]\in{\cal E}(P)$ for some $1\le \ell< p$, then
$[x_\ell,x_n]\in{\cal E}(P)$ implies
$[x_p,x_{n-1}]\in{\cal E}(P)$ by the preceding.
Hence, we may assume that $p$ is the least positive integer such that
$[x_p,x_n]\in{\cal E}(P)$.
Then Lemma~\ref{edges}(\ref{4c}) yields that any facet of $P$ that contains 
$\{x_p,x_n\}$ also contains vertices between $x_p$ and $x_n$.
Let $x_k$ be the greatest of these vertices, and let $F\in{\cal F}(P)$ such
that $\{x_p,x_k,x_n\}\subset F$.
We claim that $[x_p,x_k]\in{\cal E}(P)$ and $x_k=x_{n-1}$.

We note that (see (3.1)) $(x_k,x_n)=(y_{m-1},y_m)$,
$x_p\in\{y_0,y_{m-d+2}\}$ by Lemma~\ref{edges}(\ref{4c}), and
$[x_p,x_k]\in{\cal E}(P)$ by Lemma~\ref{edges}(\ref{4a}) and (\ref{4d}).
Since
$$E_{m-2}=[y_0,y_{m-d+1},\ldots,y_{m-3},x_k,x_n],$$
it is clear that if $x_k\ne x_{n-1}$, then $y_{m-2}=x_{k-1}$, and
there is an $F'\in{\cal F}(P)$ such that $E_{m-2}=F\cap F'$ and
$\{x_p,x_{k+1},x_n\}\subset F'$, a contradiction.
Hence $x_k=x_{n-1}$.
\end{proof}

\begin{corollary}\label{Vn}
${\cal V}_1(P)=\{x_2,\ldots, x_r\}$ and ${\cal V}_n(P)=\{x_s\ldots,x_{n-1}\}$
for some $d\le r\le n$ and $1\le s\le n-d+1$.
\end{corollary}

We are now in the position to determine some of the subpolytopes of~$P$.

\pagebreak

\begin{lemma}\label{Sd-F3}
Let $1\le s\le n-d+1$, ${\cal V}_n(P)=\{x_s,x_{s+1},\ldots, x_{n-1}\}$
and $S_d\subset\{x_s,\ldots, x_{n-1},x_n\}$.
Then
\begin{enumerate}
\item $[S_d]\subset {\cal F}(P)$ and \label{Sd}
\item $[x_0,x_{s-1},x_s,x_{n-1},x_n]\in{\cal F}_3(P)$. \label{F3}
\end{enumerate}
\end{lemma}
\begin{proof} (\ref{Sd}) Since there is an $F'\in{\cal F}(P)$ such that 
$x_n\in F'$ and $x_0\not\in F'$, it follows that $F'=[S'_d]$ for some
$S'_d\subset\{x_s,\ldots, x_n\}$ by Lemma~\ref{simplex} and the Gale property.
Since $s=n-d+1$ implies that $S'_d=\{x_s,\ldots, x_n\}$, we may assume that
$s\le n-d$.

Let $S_d\subset\{x_s,\ldots, x_{n-1}\}$ such that $|S_d\cap S'_d|=d-1$.
Then $G=[S_d\cap S'_d]\in{\cal F}_{d-2}(P)$ and $G=F\cap[S'_d]$ for some
$F\in{\cal F}(P)$.
By the Gale property, $S_d\subset F$.
Let $x_t\in S_d\subset\{x_t,\ldots, x_n\}$.
Since $[x_t,x_j]\in{\cal E}(P)$ for each $x_j\in S_d\setminus\{x_t\}$, it 
follows from Lemma~\ref{edges} that $x_t$ is the initial vertex of $F$.
Thus, $x_0\not\in F$ and $F=[S_d]$ by Lemma~\ref{simplex}.

If $s=n-d$ then we are done, and if $s<n-d$ then iterations of the preceding
argument yield (\ref{Sd}).

(\ref{F3}) Let $S_{d-4}\subseteq\{x_{s+2},\ldots,x_{n-2}\}$ and 
$S_d=\{x_s,x_{s+1}\}\cup S_{d-4}\cup\{x_{n-1},x_n\}$.
By (\ref{Sd}) and the Gale property,
$$[x_s,S_{d-4},x_{n-1},x_n]=[S_d]\cap F$$
for some $F\in{\cal F}(P)$ with $x_{s-1}\in F$.
We note that $[x_{s-1},x_n]\not\in{\cal E}(P)$ implies that $F$ is not a 
simplex, and Lemma~\ref{simplex} implies that $x_0\in F$.
Hence, $x_0=y_0$, $(x_{n-1},x_n)=(y_{m-1},y_m)$ and 
$(x_{s-1},x_s)=(y_{m-d+1},y_{m-d+2})$ by Lemma~\ref{edges}(\ref{4c}).
Finally, Lemma~\ref{edges}(\ref{4e}) with $t=m-d+1$ yields the assertion.
\end{proof}

\begin{lemma}\label{neighbors}
Let $2\le s\le n-d+1$ and ${\cal V}_n=\{x_s,\ldots, x_{n-1}\}$.
Then ${\cal V}_{n-1}(P)=\{x_{s-1},x_s,\ldots, x_{n-2},x_n\}$.
\end{lemma}
\begin{proof}
By Lemmas~\ref{between} and~\ref{Sd-F3}, 
$\{x_s,\ldots,x_{n-2},x_n\}\subset{\cal V}_{n-1}$ and
$[x_0,x_{s-1},x_s,x_{n-1},x_n]\in{\cal F}_3(P)$.
The latter and $[x_{s-1},x_n]\not\in{\cal E}(P)$ yield that 
$[x_{s-1},x_{n-1}]\in{\cal E}(P)$.
Thus, it remains to show that $[x_{s-2},x_{n-1}]\not\in{\cal E}(P)$ for 
$s\ge 3$.
We note that $[x_{s-1},x_{n-1}]\in{\cal E}(P)$ yields that there is an 
$F\in{\cal F}(P)$ such that $[x_{s-1},x_{n-1}]\subset F$, $x_s\not\in F$
and $x_{s-2}\in F$.

If $x_n\in F$ then $[x_{s-1},x_n]\not\in{\cal E}(P)$ implies that $F$ is not
a simplex and $y_0=x_0<x_{s-2}$.
We note that $(y_{m-1},y_m)=(x_{n-1},x_n)$.
Thus $x_{s-1}\in{\cal V}_{n-1}\setminus{\cal V}_n$ and 
Lemma~\ref{edges}(\ref{4c},\ref{4d}) yield that $y_{m-d}<x_{s-1}<y_{m-d+2}$.
Then $x_{s-1}=y_{m-d+1}$, $x_{s-2}=y_{m-d}$ and 
$[x_{s-2},x_{n-1}]=[y_{m-d},y_{m-1}]\not\in{\cal E}(P)$ by 
Lemma~\ref{edges}(\ref{4d}).

We may now assume that no facet of $P$ contains 
$\{x_{s-2},x_{s-1},x_{n-1},x_n\}$.
Clearly, this assumption and the Gale property yield that 
$[x_{s-2},x_{s-1},x_{n-1}]\not\in{\cal F}_2(P)$.
Thus, $[x_{s-2},x_{s-1},x_{n-1}]\subset F$ implies that $x_{n-1}=y_m$ and
$x_{s-2}\not\in T_{m-d+2}=[y_{m-d+2},\ldots, y_m]$.
Finally, $[x_{s-2},x_{n-1}]=[y_t,y_m]$ for some $0<t<m-d+2$ and 
Lemma~\ref{edges}(\ref{4c}) yield that $[x_{s-2},x_{n-1}]\not\in{\cal E}(P)$.
\end{proof}

\begin{theorem}\label{theoremb}
Let $P$ be a Gale and braxial $d$-polytope with $x_0<x_1<\cdots<x_n$,
$n\ge d+1\ge 5$ and $d$ even.
Then $P'=[x_0,x_1,\ldots,x_{n-1}]$ is a Gale and braxial $d$-polytope with
$x_0<x_1<\cdots<x_{n-1}$.
\end{theorem}
\begin{proof}
Let $F'\in{\cal F}(P')$.
We need to show that $F'$ is a $(d-1)$-braxtope and that ${\cal V}(F')$ is a 
Gale set.
We note that either $F'\in{\cal F}(P)$, or $[F',x_n]\in{\cal F}(P)$, or
$F'\not\in{\cal F}(P)$ and $x_n$ is beyond $F'$.
If $F'\in{\cal F}(P)$ then $x_n\not\in F'$ and 
$F'$ is braxial and ${\cal V}(F')$ is a Gale subset of ${\cal V}(P')$.
If $F'\not\in{\cal F}(P)$ and $[F',x_n]\in{\cal F}(P)$ then ${\cal V}(F')$ is
necessarily a Gale set and $F'$ is a $(d-1)$-braxtope by Lemma~\ref{edges}.

Let $F'\not\in{\cal F}(P)$ and $x_n$ be beyond $F'$.
We recall that ${\cal V}_n=\{x_s,\ldots, x_{n-1}\}$ for some $1\le s\le n-d+1$,
and note that ${\cal V}(F')\subset {\cal V}(P)$ and 
Lemma~\ref{inductive}(\ref{1c}) yield that 
${\cal V}(F')\subset{\cal V}_n\cup\{x_0\}$.
Thus, if $s=n-d+1$ then $F'=[x_0,x_{n-d+1},\ldots, x_{n-1}]$, whence $F'$ is a
$(d-1)$-simplex and ${\cal V}(F')$ is a Gale set.

Let $s\le n-d$ and $Y={\cal V}(F')\cap \{x_s,\ldots, x_{n-2}\}$.
We note that $|Y|\ge d-2$, and claim that $Y$ is a paired set.
We suppose otherwise and seek a contradiction.
Since $[S_d]\in{\cal F}(P)$ for each $S_d\subset\{x_s,\ldots, x_{n-2}\}$ from 
Lemma~\ref{Sd-F3}, it follows
that for some $t$ such that $1\le t\le d/2$ there is a maximal paired subset
$S_{d-2t}$ of $Y$ and a $t$-element subset $X_t$ of $Y\setminus S_{d-2t}$ such
that no two vertices of $X_t$ are successive.

Since $S_{d-2t}\cup X_t\subset\{x_s,\ldots, x_{n-2}\}$ and $n\ge s+d$, there is
an $S_d\subset\{x_s,\ldots, x_{n-2},x_{n-1}\}$ such that 
$S_{d-2}\cup X_t\subset S_d$.
Since $[S_d]\in{\cal F}(P)$ is a simplex and $x_n$ is beneath $[S_d]$, it 
follows from Lemma~\ref{inductive}(\ref{1b}) that 
$$G=[S_{d-2t},X_t,x_n]\in{\cal F}_{d-t}(P).$$
Let $F\in{\cal F}(P)$ such that $G\subseteq F$.
Since $S_{d-2t}\cup X_t\subset\{x_s,\ldots, x_{n-2}\}$ and ${\cal V}(F)$ is a
Gale set, we obtain that there is a paired set 
$$S\subset {\cal V}(F)\cap\{x_{s-1},x_s,\ldots, x_{n-2},x_{n-1}\}$$
such that $S_{d-2t}\cup X_t\subset S$.
We note that $|S_{d-2t}\cup X_t|=d-t$ implies that $|S|\ge d$, and thus, $F$ is
not a simplex and $|S\setminus\{x_{s-1}\}|\ge d-1$.
Since $[x_j,x_n]\in{\cal E}(P)$ for $x_j\in(S\setminus\{x_{s-1}\})\cup\{x_0\}$
and $x_0\in F$ by Lemma~\ref{simplex}, the contradiction we obtain is that
$x_n$ is not a simple vertex of $F$, see Lemma~\ref{edges}(\ref{4c}).

In summary: ${\cal V}(F')\subset{\cal V}_n\cup\{x_0\}$,
$Y={\cal V}(F')\cap({\cal V}_n\setminus\{x_{n-1}\})$ is a paired set, $|Y|\ge d-2$
and $Y$ contains a paired set of cardinality at most $d-2$ by Lemma~\ref{Sd-F3}.
Hence, $Y=S_{d-2}$ and ${\cal V}(F')=\{x_0\}\cup S_{d-2}\cup\{x_{n-1}\}$ is
Gale with respect to $x_0<x_1<\cdots<x_{n-1}$.
\end{proof}

\begin{corollaryth}
Let $2\le s\le n-d+1$ and ${\cal V}_n(P)=\{x_s,\ldots, x_{n-1}\}$.
Then ${\cal V}_{n-1}(P')=\{x_{s-1},x_s,\ldots, x_{n-2}\}$
\end{corollaryth}
\begin{proof}
In view of Corollary~\ref{Vn} and Lemma~\ref{neighbors}, it remains to show that
$x_{s-2}\not\in{\cal V}_{n-1}(P')$ for $s\ge 3$.

Let $[x_{s-2},x_{n-1}]\in{\cal E}(P')$.
Then $[x_{s-2},x_{n-1}]\not\in{\cal E}(P)$,
$[x_{s-2},x_{n-1},x_n]\not\in{\cal F}_2(P)$ and Lemma~\ref{inductive} yield
that $x_n$ is beyond each $F'\in{\cal F}(P')$ that contains $[x_{s-2},x_{n-1}]$.
Next,  $[x_{s-2},x_{n-1}]\in{\cal E}(P')$ and the Gale property of $P'$ imply
that there is an $F'\in{\cal F}(P')$ such that 
$[x_{s-2},x_{s-1},x_{n-1}]\subset F'$.
Since $[x_{s-1},x_{n-1}]\in{\cal E}(P)$ and $x_n$ is beyond $F'$, it follows
by Lemma~\ref{inductive} that $[x_{s-1},x_{n-1},x_n]\in{\cal F}_2(P)$ and
$x_{s-1}\in{\cal V}_n(P)$, a contradiction.
\end{proof}

\begin{theorem}\label{theoremc}
Let $P$ be a Gale and braxial $d$-polytope with $x_0<x_1<\cdots<x_n$,
$n\ge d+1\ge 5$ and $d$ even.
Let $P'= [x_0,x_1,\ldots, x_{n-1}]$, $F'\in{\cal F}(P')$, and
${\cal V}_n(P)=\{x_s,\ldots, x_{n-1}\}$, $2\le s\le n-d+1$.  Then
\begin{itemize}
\item $x_n\in\langle F'\rangle$ if $[x_0,x_{s-1},x_s,x_{n-1}]\subset F'$,
\item $x_n$ is beyond $F'$ if $F'\cap[x_0,x_{s-1},x_{n-1}]=[x_0,x_{n-1}]$, and
\item $x_n$ is beneath $F'$ if $[x_0,x_{s-1},x_s,x_{n-1}]\not\subset F'$ and
      $F'\cap[x_0,x_{s-1},x_{n-1}]\ne [x_0,x_{n-1}]$.
\end{itemize}
\end{theorem}
\begin{proof}
If $[x_0,,x_{s-1},x_s,x_{n-1}]\subset F'$ then Lemma~\ref{Sd-F3}(\ref{F3}) 
implies that $x_n\in\langle F'\rangle$.
Suppose $F'\cap[x_0,x_{s-1},x_{n-1}]=[x_0,x_{n-1}]$.
Then Lemma~\ref{edges}(\ref{4c}) and 
${\cal V}_{n-1}(P') =\{x_{s-1},\ldots, x_{n-2}\}$ yield that 
\begin{equation}\label{cardinality}
|{\cal V}(F') \cap\{x_s,\ldots, x_{n-1}\}|=d-1.
\end{equation}
Now Lemma~\ref{edges}(\ref{4c}) yields also that $[F',x_n]\not\in{\cal F}(P)$;
that is, $x_n\not\in\langle F'\rangle$.
Finally, $\langle F'\rangle\cap\{x_{s-1},x_n\}=\emptyset$, the Gale property
of $P$ and (\ref{cardinality}) yield that $F'\not\in{\cal F}(P)$.
Hence, $x_n$ is beyond $F'$.

Now assume 
\begin{equation}\label{eq-i}
[x_0,x_{s-1},x_s,x_{n-1}]\not\subset F'
\end{equation}
and
\begin{equation}\label{eq-ii}
F'\cap[x_0,x_{s-1},x_{n-1}]\ne[x_0,x_{n-1}].
\end{equation}
Thus, either $x_0\not\in F'$ or $F'\cap\{x_0,x_{n-1}\}=\{x_0\}$ or
$F'\cap[x_0,x_{s-1},x_s,x_{n-1}]=[x_0,x_{s-1},x_{n-1}]$.

We note that $[F',x_n]$ is not a $(d-1)$-simplex, and recall from 
Lemma~\ref{simplex}
that if a facet of $P'$ or $P$ does not contain $x_0$ then it is a simplex.

Let $x_0\not\in F'$.
Then $x_0\not\in[F',x_n]$ and it follows that $F'$ is a simplex and
$[F',x_n]\not\in{\cal F}(P)$; that is, $x_n\not\in\langle F'\rangle$.
Now if ${\cal V}(F')\subset{\cal V}_n$ then $F'\in{\cal F}(P)$ by
Lemma~\ref{Sd-F3}(\ref{Sd}), and if ${\cal V}(F')\not\subset{\cal V}_n$ then
$F'\in{\cal F}(P)$ by Lemma~\ref{inductive}(\ref{1c} and \ref{1a}),
so $x_n$ is beneath $F'$.

Let $x_0\in F'$.
We suppose that $x_n\in\langle F'\rangle$, and seek a contradiction.

Since $F'\in{\cal F}(P')$, there is a greatest vertex $x_p$ of $F'$ and 
$x_p\le x_{n-1}$.
If $x_p=x_{n-1}$ then $\{x_0,x_{s-1},x_{n-1},x_n\}\subset \langle F'\rangle$
by (\ref{eq-ii}), $x_s\in\langle F'\rangle$ by Lemma~\ref{Sd-F3}(\ref{F3}), and
we have a contradiction by (\ref{eq-i}).
Let $x_p<x_{n-1}$.
With reference to (\ref{ridges2}), we note that $(z_0,z_{u-1},z_u)=(x_0,x_p,x_n)$
for the facet $[F',x_n]$ of $P$
and that $z_0<z_1<\cdots<z_{u-1}$ is the vertex array of $F'$ as a facet of 
$P'$.
Next
\begin{equation}\label{eq-iii}
x_n=z_u\in\langle x_0,z_{u-d+1},z_{u-d+2},x_p\rangle
\end{equation}
by Lemma~\ref{edges}(\ref{4e}), and 
$$E'_{u-2}=[z_0,z_{u-d+1},\ldots, z_{u-3},z_{u-1}]=F'\cap F''$$
for some $F''\in{\cal F}(P')$.
We observe that $x_n\in\langle F''\rangle$ by (\ref{eq-iii}), and that 
$z_{u-1}=x_p<x_{n-1}$ and the Gale property yield that 
$x_{p-1}=z_{u-2}\not\in F''$ and $x_{p+1}\in F''$.
Clearly, $x_{p+1}$ is the greatest vertex of $F''$ as a facet of $P'$.
Thus $[E'_{u-2},x_n]$ is a $(d-2)$-face of the facet $[F'',x_n]$ of $P$
with $d=(d-1)+1$ vertices.
Since $x_{p+1}$ is not in $[E'_{u-2},x_n]$ and 
$x_0=z_0<z_{u-1}=x_p<x_{p+1}<x_n$ in the $(d-1)$-braxtope $F''$, this is
a contradiction.

In summary, $x_0\in F'$, $x_n\not\in\langle F'\rangle$ and either 
$x_{n-1}\not\in F'$ or $\{x_{s-1},x_{n-1}\}\subset F'$.
If $x_{s-1}\in F'$ then ${\cal V}(F')\not\subset{\cal V}_n\cup\{x_0\}$ and, as
already noted, $F'\in{\cal F}(P)$.
Hence, we suppose that $x_{n-1}\not\in F'$ and 
${\cal V}(F')\subset {\cal V}_n\cup\{x_0\}$.
Then ${\cal V}(F')\subset {\cal V}_{n-1}(P')\cup\{x_0\}$ and it follows from
Lemmas~\ref{V0} and~\ref{between} that any two vertices of $F'$ determine an 
edge of $F'$.
Finally, Lemma~\ref{edges} and the Gale property yield that $F'$ is a simplex
and $F'=[S_d]$ for some $S_d\subset\{x_0\}\cup\{x_s,\ldots, x_{n-2}\}$.
Since $s\ge 2$, we have a contradiction.
\end{proof}

\begin{theorem}
Let $P$ be a Gale and braxial $d$-polytope with $x_0<x_1<\cdots<x_n$,
$n\ge d+1\ge 7$ and $d$ even.
Then for some $1\le s\le n-d+1$,
$[x_j,x_n]\in{\cal E}(P)$ if and only if $j=0$ or $s\le j\le n-1$.
Furthermore
\begin{itemize}
\item $P$ is cyclic with $x_0<x_1<\cdots<x_n$ if $s=1$,
\item $P$ is periodically-cyclic with period $k=n-s+2$ if\/ $2\le s\le n-d$, 
      and
\item $P$ is a $d$-braxtope with $x_0<x_1<\cdots<x_n$ if $s=n-d+1$.
\end{itemize}
\end{theorem}
\begin{proof}
The first claim of the theorem is from Lemma~\ref{V0} and Corollary~\ref{Vn}.

Suppose $s=1$.
Then $[x_i,x_j]\in{\cal E}(P)$ for any $x_i\ne x_j$ by Lemmas~\ref{V0} 
and~\ref{between}.
From this and Lemma~\ref{edges}, it follows that each facet of $P$ is a simplex.
Since $P$ is Gale and simplicial, it is a cyclic polytope.

Let $2\le s\le n-d$ and $k=n-s+2$.
We observe that repeated applications of Theorem~\ref{theoremb} and its 
corollary yield that $\tilde{P}=[x_0,x_1,\ldots, x_{n-s+1}=x_{k-1}]$ is a Gale and
braxial $d$-polytope with $k-1\ge d+1$ and 
${\cal V}_{k-1}(\tilde{P})=\{x_1,\ldots, x_{k-2}\}$.
Thus $\tilde{P}$ is cyclic, and it readily follows from Proposition~\ref{cyclic}
and Theorem~\ref{theoremc} that $P$ is periodically-cyclic with period $k$.

Finally, let $s=n-d+1$ and $P_r=[x_0,x_1,\ldots, x_r]$, $d\le r\le n$.
We note that $P_d$ is a $d$-simplex, and thus, it is a $d$-braxtope.
Next, ${\cal V}_n(P)=\{x_{n-d+1},\ldots, x_{n-1}\}$ and repeated applications
of Theorem~\ref{theoremb} and its corollary yield that 
${\cal V}_r(P_r)=\{x_{r-d+1},\ldots, x_{r-1}\}$ for $d+1\le r\le n$.
We recall that $x_r$ is beyond $F\in{\cal F}(P_{r-1})$ only if 
${\cal V}(F)\subset{\cal V}_r(P_r)\cup\{x_0\}$, and thus, $x_r$ is beyond only
$[x_0,x_{r-d+1},\ldots, x_{r-1}]\in{\cal F}(P_{r-1})$.
With the assumption that $P_{r-1}$ is a $d$-braxtope, it is now easy to check
that Theorem~\ref{theoremc} and the preceding yield that $P_r$ is a $d$-braxtope
for $d<r\le n$.
\end{proof}

In fact, we have shown that the Gale and braxial polytopes in even dimension
$d\ge 6$  are exactly
the polytopes of Proposition~\ref{cyclic}, that is, the periodically-cyclic
Gale polytopes constructed in \cite{bisz-pcg}.
Henceforth call these polytopes {\em Gale-braxial polytopes}.
In view of the comments preceding Proposition~\ref{cyclic},
we conjecture that, for even $d\ge 6$, there exist periodically-cyclic 
$d$-polytopes that are Gale, but not Gale-braxial.

\begin{problem}{\em
Determine all periodically-cyclic $d$-polytopes for even $d\ge 6$.
}\end{problem}

In \cite{bisz-pcg} it was shown that the construction of 
Proposition~\ref{cyclic}, applied in dimension four, produces polytopes that
are not periodically-cyclic.
Theorems~\ref{theorema}--\ref{theoremc} show that Gale and braxial
4-polytopes are polytopes constructed in this manner, so they are not
periodically-cyclic.
As mentioned in Section~2, some bicyclic polytopes are both Gale and
periodically-cyclic (and necessarily not braxial) \cite{bisz-pcg4}.
Naturally, we would like a better understanding of braxial polytopes
with regard to these other properties.

\begin{problem}{\em
Determine the explicit facet structure of Gale-braxial 4-polytopes.
}\end{problem}

\begin{problem}{\em Determine all 4-polytopes that are braxial and 
periodically-cyclic.}
\end{problem}

Of course, we can ask for a complete classification of braxial 4-polytopes
but this is not likely to be tractable.

We have focused here on even-dimensional polytopes.  
In odd dimensions ordinary polytopes have been studied as a generalization
of cyclic polytopes.  
The question of whether ordinary polytopes are periodically cyclic is open;
Dinh \cite{dinh} proved that ordinary polytopes satisfy a weaker condition,
local neighborliness.


\begin{thebibliography}{10}

\bibitem{bayer-bisz}
M.~M. Bayer and T.~Bisztriczky.
\newblock On braxtopes, a class of generalized simplices.
\newblock arXiv:math.CO/0607396, 2006.

\bibitem{bisz-mult}
T.~Bisztriczky.
\newblock On a class of generalized simplices.
\newblock {\em Mathematika}, 43:274--285, 1996.

\bibitem{bisz-ord}
T.~Bisztriczky.
\newblock Ordinary $(2m+1)$-polytopes.
\newblock {\em Israel J. Math.}, 102:101--123, 1997.

\bibitem{bisz-pcg}
T.~Bisztriczky.
\newblock A construction for periodically-cyclic {G}ale $2m$-polytopes.
\newblock {\em Beitr\"{a}ge Algebra Geom.}, 42(1):89--101, 2001.

\bibitem{bisz-multrig}
T.~Bisztriczky and K.~B{\"o}r{\"o}czky.
\newblock Oriented matroid rigidity of multiplices.
\newblock {\em Discrete Comput. Geom.}, 24(2-3):177--184, 2000.
\newblock The Branko Gr\"unbaum birthday issue.

\bibitem{bisz-pcg4}
T.~Bisztriczky and K.~B{\"o}r{\"o}czky.
\newblock On periodically-cyclic Gale 4-polytopes.
\newblock {\em Discrete Math.} 241:103--118, 2001.
\newblock Selected papers in honor of Helge Tverberg.

\bibitem{dinh}
T.~Dinh.
\newblock {\em Ordinary Polytopes.}
\newblock Ph. D. Thesis, The University of Calgary, 1999.

\bibitem{gale}
D.~Gale.
\newblock Neighborly and cyclic polytopes.
\newblock In {\em Proc. Sympos. Pure Math., Vol. VII}, pages 225--232. Amer.
  Math. Soc., Providence, R.I., 1963.

\bibitem{grunbaum}
B.~Gr{\"u}nbaum.
\newblock {\em Convex polytopes}, volume 221 of {\em Graduate Texts in
  Mathematics}.
\newblock Springer-Verlag, New York, second edition, 2003.
\newblock Prepared and with a preface by Volker Kaibel, Victor Klee and
  G\"unter M.\ Ziegler.

\bibitem{klee-minty}
V. Klee and G. J. Minty.
\newblock How good is the simplex algorithm?
\newblock In {\em Inequalities, III (Proc. Third Sympos., Univ. California, Los
  Angeles, Calif., 1969; dedicated to the memory of Theodore S. Motzkin)},
  pages 159--175. Academic Press, New York, 1972.

\bibitem{smilansky}
Z. Smilansky.
\newblock Bi-cyclic {$4$}-polytopes.
\newblock {\em Israel J. Math.}, 70(1):82--92, 1990.

\bibitem{ziegler}
G.~Ziegler.
\newblock {\em Lectures on polytopes}, volume 152 of {\em Graduate Texts in
  Mathematics}.
\newblock Springer-Verlag, New York, 1995.

\end{thebibliography}
\end{document}